\documentclass{amsart}
\usepackage{graphicx, enumerate}
\usepackage{verbatim} 
\usepackage{syntonly}
\usepackage{url}  

\newtheorem*{namedtheorem}{\theoremname}
\newcommand{\theoremname}{testing}

\begin{document}

\title {Fermat's Four Squares Theorem} 

\author{Alf van der Poorten}
\address{Centre for Number Theory Research, 1 Bimbil Place, Killara, Sydney,
NSW 2071, Australia}
\email{alf@maths.usyd.edu.au (Alf van der Poorten)}



\renewcommand{\abstractname}{Introduction}
\begin{abstract} It is easy to find a right-angled
triangle with integer sides whose area is $6$. There  is no such triangle with
area $5$, but there is one with rational sides (a `\emph{Pythagorean triangle}'). 

For historical
reasons, integers such as $6$ or $5$ that are (the squarefree part of) the area
of some Pythagorean triangle are called `\emph{congruent numbers}'.  These
numbers actually are interesting for the following reason: Notice the
sequence $\frac14$, $6\frac14$,
$12\frac14$. It is an arithmetic progression with common difference $6$,
consisting of squares
$(\frac12)^2$,
$(\frac52)^2$,
$(\frac72)^2$ of rational numbers. Indeed the common difference of three
rational squares in AP is a congruent number and every congruent number is the common difference of three rational squares in arithmetic progression. 

The triangle given by $9^{2}+40^{2}=41^{2}$ has area $180=5\cdot6^{2}$ and the numbers $x-5$, $x$ and $x+5$ all are rational squares if $x=11\frac{97}{144}$. Recall one obtains all Pythagorean triangles with relatively prime integer sides by taking $x=4uv$, $y=\pm(4u^{2}-v^{2})$, $z=4u^{2}+v^{2}$ where $u$ and $v$ are integers with $2u$ and $v$ relatively prime.

Fermat proved that there is no AP of
more than three squares of rationals.  
\end{abstract}

\maketitle

\noindent Several years ago (in fact, at an AMSI Summer School), after I had pointed out that three rational squares in arithmetic progression with integer common difference correspond to a Pythagorean triangle with that integer area, I found it natural to tell my audience that ``Fermat proved that there is no arithmetic progression of
more than three squares (of rationals). In other words, the pair of
\emph{diophantine equations} $a^2+c^2=2b^2$ and $b^2+d^2=2c^2$ has no solution in rationals $a$, $b$, $c$ and $d$.'' But my suggestion that they be Fermat and write me an essay on the proof fell on stony ground and, worse, the best \emph{I} could provide as a solution was to say: ``Too hard for me? I looked this up in \cite{119} and found at p.\,54 the unhelpful footnote ``Fermat\ could show by descent that one cannot have four squares in AP \dots\;. Gerry Myerson has pointed me to a reference but the argument there seems utterly soulless and I remain searching for a decent descent argument warranting report to you.'' 

I decided recently that such a proof was most readily found on a (previously) blank page of my notebook.

\renewcommand{\theoremname}{Fermat's four squares theorem}
\begin{namedtheorem} 
There are no four distinct rational squares in arithmetic progression.
\end{namedtheorem}

Fermat's four squares theorem seems to appear in the literature as a mildly surprising corollary of other somewhat obscure diophantine results, possibly because authors start by translating the suggestion that if $r$, $s$, $u$ and $t$ are squares in arithmetic progression then $t-u=u-s=s-r$, giving a pair of equations each involving three squares. 

Here I give a more direct proof, starting from four integer squares $x-6n$, $x-2n$, $x+2n$ and $x+6n$ and remarking that I may suppose without loss of generality that the four squares all are odd and hence have a common difference, here $4n$, divisible by $4$. Thus $x$ is odd and, more, we may suppose that the four squares are pairwise relatively prime. Plainly we also have an odd integer $y$ prime to $x$ and $n$ so that
$$
y^2=(x^2-4n^2)(x^2-36n^2)=(x^2-20n^2)^2-256n^4\,.
$$
It follows that we have a Pythagorean triple $(16n^2,y,x^2-20n^2)$ and hence relatively prime integers $2u$ and $v$ so that $4uv=16n^2$, $4u^2+v^2=x^2-20n^2$. Hence $u$ must be even and there are integers $A$ and $D$, with $D$ odd, so that $u=4A^2$ and $v=D^2$. We also have
$$4u^2+v^2=x^2-5uv \quad\text{so}\quad (4u+v)(u+v)=x^2$$
implying that both $4u+v=16A^2+D^2$ and $u+v=4A^2+D^2$ are squares.

Because each of the implications above is reversible we have in fact shown that there are relatively prime integers $A$ and $D$ so that $16A^2+D^2$ and $4A^2+D^2$ both are squares if and only if there are four pairwise relatively prime squares in arithmetic progression with common difference $4AD$.

However, $4A^2+D^2$ a square entails there are relatively prime integers $U$ and $V$ so that $2UV=2A$ and $U^2-V^2=\pm D$, while $16A^2+D^2$ a square yields relatively prime integers $2U'$ and $V'$ so that $4U'V'=4A$ and $4U'{}^2-V'{}^2=\pm D$.

In particular, the two different factoris{ations} $UV=U'V'$ of the even integer $A$ entail there are pairwise relatively prime integers $2a$, $b$, $c$, $d$ so that $\pm D=4a^2b^2-c^2d^2=16a^2c^2-b^2d^2$. That is, $b^2(4a^2+d^2)=c^2(16a^2+d^2)$ and one sees that both $4a^2+d^2$ and $16a^2+d^2$ are squares.

So there are four pairwise relatively prime squares in arithmetic progression with common difference $4ad$. However, $ad$ is a proper divisor of $A$ and thus is certainly smaller than $AD$, proving Fermat's four squares theorem by descent.

Here I have tacitly supposed that the given arithmetic progression is nontrivial; that is that $n\ne0$. That tacit presumption is of course important because $1$, $1$, $1$, $1$ \emph{is} an arithmetic progression of pairwise relatively prime squares; albeit a degenerate such progression. However, I use the tacit assumption in an important manner only in my final paragraphs. If $n=0$ then necessarily $A=0$ and $D^2=1$ and those paragraphs, corrected for the degenerate case, do not descend but appropriately reproduce the given trivial progression.

\subsection*{What's going on here?} One can rewrite the opening assumption as alleging that the curve $\mathcal C:Y^{2}-(X^{2}-5)Y+4=0$ contains a rational point $(X,Y)$; specifically, so that $X$ has denominator $2n$ and $Y$ has denominator $4n^{2}$.

Indeed, $\mathcal C$ is a quartic model for an elliptic curve $\mathcal E: y^{2}=x(x+1)(x+4)$ obtained by taking $x=Y$ and $y=XY$; thus, for the presumed rational point on $\mathcal E$, the denominator of $x$ is $4n^{2}$ and that of $y$ is $8n^{3}$. $\mathcal E$ is curve 24A1 of Cremona's tables \cite{Cr}. My argument  confirms that there are no rational points on $\mathcal E$ corresponding to a nontrivial arithmetic progression.

\subsection*{History}  After concocting the remarks above, I checked Dickson{'s} \textit{History of the Theory of Numbers}~\cite{Di}. Dickson reports at II, XIV, p.\,440 that Fermat proposed the problem of constructing a nontrivial sequence of four squares of rationals in arithmetic progression to Frenicle in 1640 and stated that it is impossible; and \textit{inter alia} gives a summary of an uncompelling 1813 proof. Dickson, at II, XXII, p.\,635 mentions an argument of Euler  which leads one to see that $x^{2}+y^{2}$ and $x^{2}+4y^{2}$ are not both squares for $x$ odd, $y\ne0$ even; with the four squares theorem a corollary. I imagine that this effectively coincides with my argument. Reassuringly, back at p.\,440 Dickson cites an 1898 \textit{Amer.\ Math.\ Monthly} item --- which, at \textbf 5, p.\,180, turns out to be the problem ``Find, if possible, four square numbers in arithmetical progression'' --- and drily remarks that ``Several writers failed to find a solution.''

\bibliographystyle{amsalpha}

\label{page:lastpage}
\end{document}